\documentclass[12pt]{article}
\usepackage{amsmath}
\usepackage{graphicx}
\usepackage{comment} 
\usepackage{xcolor}
\usepackage{amssymb}
\usepackage{amsthm}

\newcommand{\blind}{0}

\begin{document}

\def\spacingset#1{\renewcommand{\baselinestretch}%
{#1}\small\normalsize} \spacingset{1}


\if0\blind

  \title{\bf Noise-induced synchronization and regularity in feed-forward-loop motifs}
 \author{Gurpreet Jagdev$^a$,  Na Yu$^a$\\
    {\small $^a$Department of Mathematics, Toronto Metropolitan University, Toronto, Canada}\\ }
  
  \date{}    
  \maketitle
 \fi

\if1\blind
{
  \bigskip
  \bigskip
  \bigskip
  \begin{center}
    {\LARGE\bf Title}
\end{center}
  \medskip
} \fi

\bigskip


\spacingset{1.8} 
\section{Introduction}

Network motifs, a fundamental concept in network analysis, are recurring interaction patterns observed in various systems \cite{Mangan2003}. These motifs often retain their specific dynamical functions as well as integrate into the collective behaviors of the broader network structures, thus they are often considered as the building blocks of complex networks \cite{Shen2002, reigl2004search}. Among the common three-node network motifs, coherent and incoherent feedforward loops (FFL) have proven to be notably more prevalent than their counterparts  \cite{reigl2004search, Song2005}. Here we call them type-1 motif (T1) and type-2 motif (T2) for simplicity and they are demonstrated in Fig. 1.  In motif T1, three unidirectional connections are positive (i.e., excitatory coupling). In motif T2, the connections from node 1 to node 2 and from node 1 to node 3 are positive but the connection from node 2 to node 3 is negative (i.e., inhibitory coupling). 

FFL network motifs have been found in many actual biological networks, such as gene expression in bacteria and yeast \cite{Milo2002, Lee2002}, human and mouse genomes \cite{Boyer2005, tsang2007}, cat cortex \cite{sporns2004}, and the nervous system of roundworm \cite{milo2002network}. Moreover, both motifs are constructed such that node 1 may be considered an input layer, and node 3 an output layer of the network. There are two parallel signalling pathways: one direct pathway from the input layer to the output layer, and one indirect pathway from input to output by means of node 2 (the layer of internodes or interneurons). Such parallel information transmission structure has been discovered in the auditory cortex \cite{Eggermont1998} and electrosensory system \cite{Middleton2006}.

The significance of FFL network motifs leads to the question of how they interact with one another and perform specific information-processing roles in the presence of noise, a ubiquitous factor in real life. Under the conditions of equal coupling and symmetric noise, FFL motifs have been shown to produce coherence \cite{Gui2016} or resonance \cite{Krauss2019}. However, to our best knowledge, no prior have clearly demonstrated their performance in a heterogeneous setting: unequal coupling and asymmetric noise. To address this unexplored aspect, our research will focus on noise-induced synchrony and regularity of T1 and T2 motifs in the heterogeneous setting. The remainder of this paper is structured as follows. Sec. 2 introduces the mathematical model and methods. Sec. 3.1 discusses the noise-induced dynamics of our networks. Sec. 3.2 studies the effects of noise on network synchronization. Sec. 3.3 studies the effects of noise on output regularity. Sec. 3.4 studies the effects of network connectivity on network synchrony and output regularity. Sec. 3.5 studies the effects of the bifurcation parameter $\lambda_0$ on network synchronization and output regularity. And a discussion is given in Sec. 4.

\begin{figure}
    \centering
    \includegraphics[scale=0.8]{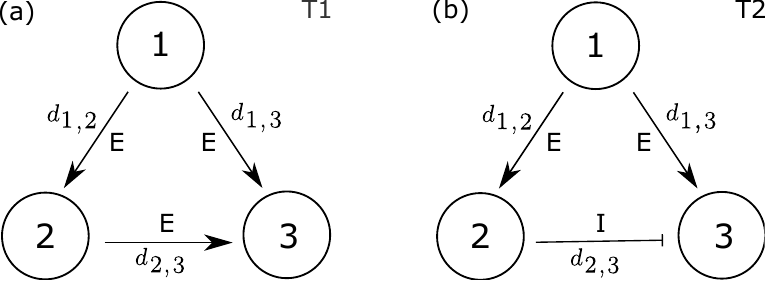}
    \caption{Coherent and incoherent feedforward-loop (FFL) motifs: (a) motif type 1 (T1); and (b) motif type 2 (T2). The term $d_{j,i}$ denotes the coupling strength, where $d_{1,2}>0, d_{1,3}>0$ in T1 and T2, $d_{1,3}>0$ in T1, and $d_{1,3}<0$ in T2. The letters ``E" and ``I" denote excitatory and inhibitory coupling.}
    \label{fig:motif}
\end{figure}

\section{Mathematical model and methods}
\subsection{Model}

We consider a system comprised of three coupled oscillators arranged in FFL patterns (see Fig. 1), where the dynamics of each oscillator are described by the canonical model for the normal form of a Hopf bifurcation (HB), the $\lambda-\omega$ system, with additive noise and diffusive coupling terms. The $i$th oscillator is modelled by the system of stochastic differential equations (SDEs)
\begin{align}
dx_i &= \Big[\lambda(r_i)x_i - \omega(r_i)y_i + \sum_{j \neq i}d_{j,i}(x_j-x_i)\Big]dt+\delta_i d\eta_i(t), \label{x_i}\\
dy_i & = \Big[\omega(r_i)x_i + \lambda(r_i)y_i + \sum_{j \neq i} d_{j,i}(y_j-y_i)\Big]dt, \label{y_i}
\end{align}
where $i,j=1,2,3$. 

The amplitude and phase of $i$th oscillator can be obtained by $r_i = \sqrt{x_i^2 + y_i^2}$ and $\phi_i = \tan^{-1} \frac{y_i}{x_i}$, respectively. The function $\lambda(r_i) = \lambda_0 + \alpha r_i^2 + \gamma r_i^4$ governs the modulation of the amplitude of $i$th oscillator. $\lambda_0$ serves as the control parameter, and the Hopf bifurcation (HB) occurs when $\lambda_0=0$. The parameters $\alpha$ and $\gamma$ influence the system's behavior away from this bifurcation point. Additionally, the function $\omega(r_i) = \omega_0 + \omega_1 r_i^2$ controls the modulation of the $i$th oscillator's frequency, with $\omega_1$ governing how the frequency evolves concerning the amplitude $r_i$. Notably, when $\omega_1=0$, changes in amplitude do not directly impact the phase. We consider a supercritical Hopf bifurcation and the corresponding parameter values are $\alpha=-0.2$, $\gamma=-0.2$, $\omega_0=2$, and $\omega_1=0$. The term $\delta_i d\eta_i(t)$ represents intrinsic white noise applied to $x_i$, and $\eta_i(t)$ characterizes a Wiener process with zero mean and unit variance. The parameter $\delta_i$ is the noise intensity.

Additionally, the terms $\sum_{j \neq i}d_{i,j}(x_j-x_i)$ and $\sum_{j \neq i} d_{i,j}(y_j-y_i)$ denote the diffusive coupling that one neuronal oscillator receives, where $d_{i,j}$ signifies the connection strength between oscillators $i$ and $j$, indicating the signal transmitted from oscillator $i$ to oscillator $j$. To create the Feed-Forward Loop (FFL) motifs, as illustrated in Figure 1, we configure $d_{j,1}=0$ and $d_{3,i}=0$ for all $i$ and $j$. This setup ensures that oscillator 1 does not receive input from other oscillators within the FFL, and oscillator 3 does not transmit output to other oscillators within the FFL. In motif T1, all three connections are excitatory, meaning that $d_{1,2}>0$, $d_{1,3}>0$, and $d_{2,3}>0$. Conversely, in motif T2, there are two excitatory connections ($d_{1,2}>0$ and $d_{1,3}>0$) along with one inhibitory connection ($d_{2,3}<0$).

\subsection{Methods}

All computational tasks, including simulation, numerical analysis, and figure generation, are carried out using MATLAB. We employ the Euler-Maruyama method to calculate the numerical solutions to SDEs and use the time interval $[50, 200]$ and time step of $dt=0.01$. We initiate the simulations with arbitrary, small random initial conditions $x_i(0),y_i(0) \sim N(0,0.008^2)$ for $i=1,2,3$. To address the challenges posed by the high-frequency fluctuations inherent in noise-induced oscillations during numerical analysis, we implement a low-pass filter. This filter is constructed using a Gaussian-weighted moving average spanning a window of $100$ data points, ensuring enhanced result consistency. Finally, we compute the coherence measures $\sigma$, $\gamma$, and $R$ outlined in Equations 4 through 6, averaging these values over a total of $N=200$ simulations.


\section{Results}
 \subsection{Noise-induced dynamics}
 \label{sec: noise-induced dynamics}
\begin{figure}
    \centering
    \includegraphics[scale=0.81]{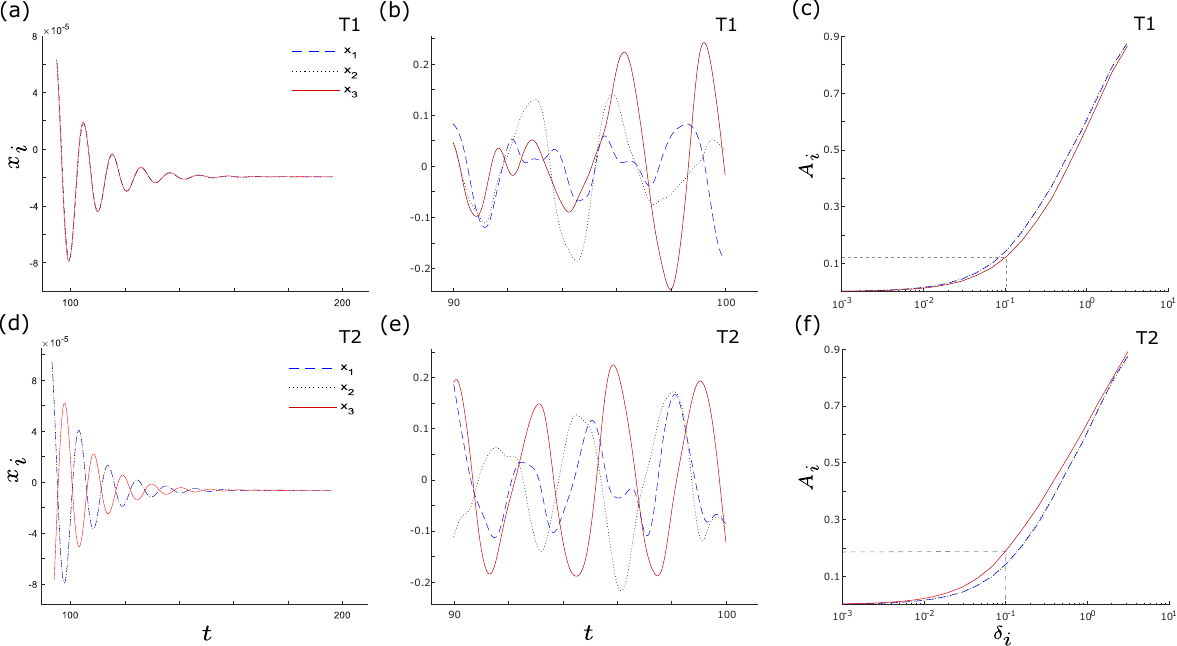}
    \caption{(a) and (d): timeseries of $x_i$, for $\delta_i=0$, $i=1,2,3$, for motifs $T1$ and $T2$, respectively. (b) and (e): timeseries of $x_i$ for $\delta_i=0.01$, $i=1,2,3$, for motifs $T1$ and $T2$, respectively. (c) and (f): amplitude, $A_i$, of $x_i$ as a function of $\delta_i$, for motifs $T1$ and $T2$, respectively. Dashed blue lines represent $x_1$, dotted black lines $x_2$, and solid red lines $x_3$. Other parameters are: $\alpha=-0.2$; $\gamma=-0.2$; $\omega_0=2$; $\omega_1=0$; $\lambda_0=-0.3$; $d_{3,1}=d_{3,2}=d_{2,1}=0.01$ for T1; and $d_{3,1}=-d_{3,2}=d_{2,1}=0.01$ for T2.}
    \label{fig:fig2}
\end{figure}

We consider the dynamics of our network motifs in the excitable regime but in the vicinity of a supercritical HB at $\lambda_0=0$ (e.g., $\lambda_0=-0.3$ in Fig. 2). The deterministic systems associated with motifs T1 and T2 ($\delta_i=0$ for $i=1,2,3$) converge to the stable fixed point $(0,0)$, as shown in Fig. 2a and 2d, respectively. We also observed that all three oscillators of the T1 motif reach in-phase during the transient time, despite their different initial values (Fig. 2a). The deterministic T2 oscillators 1 and 2 exhibit in-phase, whereas oscillator 3 show anti-phase with the others during the transient time (Fig. 2d).

The presence of intrinsic noise induces sustained limit-cycle oscillations in all oscillators within the T1 and T2 motifs (Fig. 2b and 2e where $\delta_i$ equals 0.01 for $i=1,2,3$). This kind of oscillation is commonly referred to as "noise-induced oscillation". Examining the dynamics of the two input oscillators (specifically, oscillator 1) in both T1 and T2 motifs, we observe similarities in their periods and time-varying amplitudes, as indicated by the blue dashed lines in Fig. 2b and 2e. A similar resemblance is found for oscillator 2, as depicted by the gray dotted lines. However, a notable distinction emerges when we consider the output oscillators (oscillator 3) of the T1 and T2 motifs.  The output oscillator of the T2 motif (red curve in Fig. 2e) exhibits more regular oscillations compared to its counterpart in the T1 motif (red curve in Fig. 2b). For example, the amplitude of the output oscillator in motif T2 remains relatively stable, characterized by more pronounced peaks and fewer small-amplitude perturbations.

To determine whether this is the case for other noise intensities, we compute the time-averaged amplitude of $x_i$, $A_i = \frac{1}{T}\int^T_0|x_i(t)|\ dt,$ $i=1,2,3$, as a function of $\delta_i$, and average over $N=200$ trials. Our results are displayed in Fig. 2c for motif T1 and Fig. 2f for motif T2. The dashed blue lines represent $A_1$, the dotted black lines represent $A_2$, and the solid red lines represent $A_3$. We find that $A_1$ and $A_2$ are overlapped for both motifs T1 and T2, implying that oscillators 1 and 2 

However, $A_3$ differs between the T1 and T2 motifs. In particular, for motif T2, $A_3$ is greater than its counterpart in motif T1. For instance, when $\delta_1,\delta_2,\delta_3=0.1$, $A_3 \approx 0.12$ for the T1 motif and $A_3 \approx 0.19$ for the T2 motif (see the dashed grey lines in Figs. 2c and 2f). Furthermore, not only does $A_3$ itself differ across the two FFL motifs, but the relative size of $A_3$ compared to $A_1$ and $A_2$ differs across the two FFL motifs as well. In particular, for motif T2 $A_3 \geq A_1$, $A_2$, whereas $A_3 \leq A_1$, $A_2$ for motif T1.

\subsection{The effects of noise on network synchronization}
\label{sec:network sync}
Next we consider the impact of noise on the synchronization of noise-induced dynamics. To excite oscillations from rest, we fix the noise intensity applied to $x_2$ and $x_3$ (e.g., $\delta_2,\delta_3=0.01$ in Fig. 3) and vary the noise intensity applied to $x_1$ (e.g. $\delta_3 \in [0.001,10^{1/2}]$ in Fig. 3). Then, we consider network synchronization from two perspectives. We begin first from a broad perspective, where we analyze the synchronization between every oscillator in each FFL motif. In order to quantify the degree of network synchrony, we use the root mean square deviation \cite{wang2009,gao2001}
\begin{equation}
\label{sigma}
    \sigma = \frac{1}{T}\int^T_{0}\sigma_t\ dt,
\end{equation}
where
\begin{equation}
\label{sigma(t)}
  \sigma_t = \sqrt{\frac{1}{M}\sum_{i=1}^M \left( \frac{x_i(t)}{A_i}\right)^2 - \left(\frac{1}{M}\sum_{i=1}^M \frac{x_i
  (t)}{A_i}\right)^2}
\end{equation}
and $M$ is the total number of oscillators. Note that $\sigma_t$ is computed using the amplitude-normalized time series, $x_i(t)/A_i$, rather than $x_i(t)$. This normalization ensures that $\sigma$ is a measure of temporal synchronization, or phase synchronization, rather than complete synchronization, which depends on alignments in both amplitude and phase \cite{rosenblum2001phase}. Furthermore, since $\sigma$ quantifies the degree of variability between $x_1$, $x_2$, and $x_3$, it follows that smaller values of $\sigma$ indicate greater levels of synchrony.

\begin{figure}
    \centering
    \includegraphics[scale=0.7]{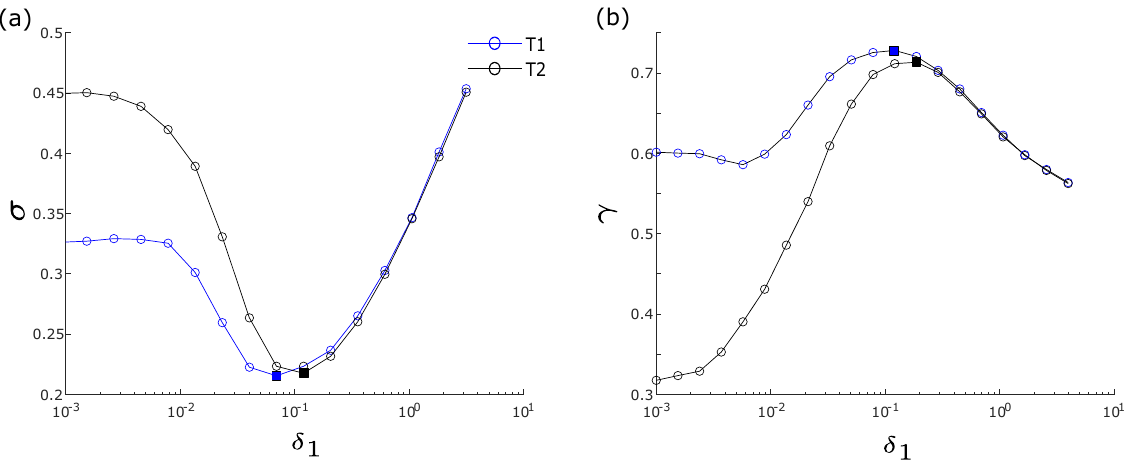}
    \caption{(a) root mean square deviation, $\sigma$, vs. noise intensity, $\delta_1$. (b) mean phase coherence, $\gamma$, vs. noise intensity, $\delta_1$. Blue lines represent the T1 motif, black lines represent the T2 motif, and solid squares indicate the minimum and maximum values of $\sigma$ and $\gamma$, respectively. The T1 motif has $d_{3,1}=d_{3,2}=d_{2,1}=0.1$, and the T2 motif has $d_{3,1}=-d_{3,2}=d_{2,1}=0.1$. Other parameters are $\alpha=-0.2$, $\gamma=-0.2$, $\omega_0=2$, $\omega_1=0$, $\delta_2=\delta_3=0.01$, and $\lambda_0=-0.1$.}
    \label{fig:network_sync}
\end{figure}

The results of our analysis are presented in Fig.  3a which displays the root mean square deviation, $\sigma$, as a function of the driving noise intensity, $\delta_1$, for the T1 (blue curve) and T2 (black curve) motifs. We observe that both of the T1 and T2 motifs display a resonant response to $\delta_1$, which indicates the occurrence of coherence resonance (CR) \cite{gang1993stochastic,pikovsky1997coherence}. That is, when $\delta_1$ is weak (e.g. $\delta_1<0.12$ for motif T2 in Fig.  3a), network synchrony increases with the increment of $\delta_1$, or equivalently, the value of $\sigma$ decreases. Then, at an intermediate intensity of $\delta_1$, for example, $\delta_1 \approx 0.12$ for motif T2 in Fig.  3a, network synchronization reaches an optimal state as $\sigma$ attains its absolute minimum. Then, as the noise intensity increases further, network synchrony exhibits a sharp decline which suggests that the noise has begun to dominate network dynamics. Indeed, our results in Fig. 3 indicate that the network synchronization (i.e. $\sigma$) is the same for the T1 and T2 motifs across the strong intensity range of $\delta_1$ (e.g. $\delta_1>0.5$ in Fig. 3).

We observe that both motifs display a resonant character and behave similarly over the strong intensity range of $\delta_1$, however, they differ considerably over the weak and intermediate intensity ranges (e.g. $\delta_1 \leq 0.1$ in Fig.  3a). Over the weak intensity range of $\delta_1$, the $\sigma$ curve which corresponds to the T1 motif (blue curve in Fig.  3a) is always less than the black curve in Fig.  3a which corresponds to the T2 motif. This suggests that the T1 motif displays comparatively greater network synchrony. Moreover, we see that the optimal driving noise intensity, denoted by $\delta_1^*$, and the absolute minimum which it induces, denoted $\sigma^*$, are both smaller for the T1 motif than the T2 motif. For example, Fig.  3a shows that the optimal driving noise intensity and its corresponding minimum are located at $\delta_1^*=0.07$ and $\sigma^*=0.215$ for motif T1, and $\delta_1^*=0.12$ and $\sigma^*=0.218$ for motif T2 (see solid squares in Fig.  3a). 

A second perspective from which to study network synchronization comes by examining the synchrony solely between the input and output oscillators of the networks. To quantify the synchronization between the input-output pair $x_1$ and $x_3$ we analyze the distribution of their phase differences, $\Delta\phi=\phi_1-\phi_3$, using a measure known as the mean phase coherence \cite{rosenblum2001phase,mormann2000mean}:
\begin{equation}
\label{gamma}
\gamma= \sqrt{\left(\frac{1}{T} \int^{T}_{t_0} \sin{\Delta \phi} \ dt\right)^2 + \left(\frac{1}{T} \int^{T}_{t_0} \cos{\Delta \phi} \ dt\right)^2}.
\end{equation} The range of $\gamma$ is $0 \leq \gamma \leq 1$, such that larger values of $\gamma$ indicate a greater degree of phase synchronization; in particular, $\gamma=1$ indicates a state of perfect synchronization and $\gamma=0$ indicates a completely chaotic state. 
The results of our analysis are presented in Fig.  3b, where the blue and black curves correspond to the T1 and T2 motifs, respectively. As in Fig. 2a, $\gamma$ is computed over the range $0.001 \leq \delta_1 \leq 10^{1/2}$ with $\delta_2$ and $\delta_3$ fixed at 0.01. The results of Fig.  3b corroborate those in Fig.  3a, and reveal a resonant behaviour in the $\gamma$ curves which mirrors that in Fig.  3a. That is, for weak noise intensities (e.g. $\delta_1<0.12$ for motif T1), increasing $\delta_1$ leads to an increase in network synchronization as indicated by an increase in $\gamma$. Then, at an intermediate noise level (e.g. $\delta_1= 0.12$ for motif T1), synchronization reaches an optimal state, where $\gamma$ reaches its absolute maximum. For stronger noise intensities, network synchronization declines which is indicated by a falling $\gamma$. This suggests that optimal phase synchrony between $x_1$ and $x_3$ is achieved at an intermediate level of $\delta_1$. 

Moreover, we find that the T1 motif requires a lesser intensity of noise to attain an optimal level of synchrony than the T2 motif. For example, in Fig. 3b, the optimal noise intensity for the T1 motif, $\delta_1^*=0.12$, is smaller than that for the T2 motif which is $\delta_1^*=0.188$, with a larger corresponding maximum in $\gamma$, $\gamma^*=0.728$, relative to the T2 motif, for which $\gamma^*=0.713$. Furthermore, we note that $\gamma$ for the T1 motif is consistently higher than that of the T2 motif, and that this difference diminishes as $\delta_1$ increases (as in Fig.  3a). Although the results shown in Figs. 3a and 3b have the same qualitative characteristics, there are quantitative differences such as the optimal driving noise intensities. These differences emerge naturally as a result of the inherent differences in the measures $\gamma$ and $\sigma$ defined in Eqs. \ref{gamma} and \ref{sigma}, respectively. Nevertheless, both measures indicate that the T1 motif exhibits a greater propensity for network synchronization than the T2 motif over the weak to intermediate noise intensity ranges and requires less noise to achieve an optimal level of network synchronization.

 \subsection{The effects of noise on output regularity}
 \label{sec: regularity}
 In addition to studying the effects of noise on network synchronization, we also consider the effects of noise on output regularity by analyzing the regularity of the time series $x_3$. By regularity, we refer to the degree to which the dynamics of $x_3$ are periodic and quantify this using the coefficient of variation of the inter-spike intervals (ISIs). The coefficient of variation, $R$, is defined as the standard deviation of the ISIs divided by the mean of the ISIs \cite{bonsel2022control,lu2019phase}, that is,
 \begin{equation}
\label{coeff variation}
    R = \frac{\sqrt{\frac{1}{K-1}\sum_{k=1}^{K-1} \left(\text{ISI}_k\right)^2 - \left(\frac{1}{K-1}\sum_{k=1}^{K-1}\text{ISI}_k \right)^2}}{\frac{1}{K-1}\sum_{k=1}^{K-1} \left(\text{ISI}_k\right)},
\end{equation}
where $\text{ISI}_k$ denotes the $k$th ISI. In order to compute an ISI we mark the occurrence of a spike by a peak in the time series of $x_3$ (e.g., inset of Fig.  4a) and compute the difference between consecutive spike times. Furthermore, because $R$ is a measure of central variability, larger values of $R$ indicate a more irregular firing pattern and consequently lower regularity, and vice versa. 
 
 \begin{figure}
     \centering
    \includegraphics[scale=0.75]{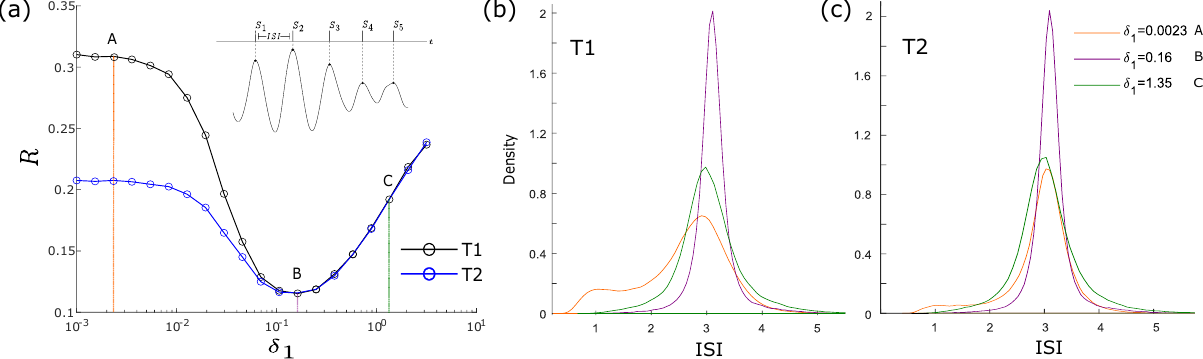}
    \caption{(a): coefficient of variation, $R$,  vs. noise intensity, $\delta_1$, for motifs T1 (black curve) and T2 (blue curve). Points A, B, and C correspond to noise intensities of $\delta_1 =$ $0.0023$, $0.16$, and $1.35$. Inset of panel (a): inter-spike interval calculation. (b) and (c): density functions of inter-spike intervals for $\delta_1$ corresponding to points A (orange line), B (purple lines), and C (green line) for the T1 and T2 motifs, respectively. Other parameters are $\alpha=-0.2$, $\gamma=-0.2$, $\omega_0=2$, $\omega_1=0$, $\delta_2=\delta_3=0.01$, and $\lambda_0=-0.1$, with $d_{31}=d_{32}=d_{21}=0.1$ for the T1 motif and $d_{31}=-d_{32}=d_{21}=0.1$ for the T2 motif.}
    \label{fig:regularity}
 \end{figure}
We present the results of our analysis in Fig.  4a, where $R$ is plotted against the driving noise intensity, $\delta_1$, for the T1 motif (black curve) and T2 motif (blue curve). We find the that both motifs display a resonant character, which indicates the occurrence of CR. That is, $R$ is relatively low at weak intensities of $\delta_1$ (e.g. $\delta_1=0.0023$ in Fig. 4a), and with an increase in $\delta_1$, $R$ decreases, indicating that output regularity is increasing, or, the firing pattern is becoming more regular. Then, as $\delta_1$ surpasses an optimal intensity, $R$ starts to increase again, and the output of our networks becomes more chaotic. Thus, the addition of intrinsic noise can optimize the regularity of the output oscillator, $x_3$, of both the T1 and T2 motifs at an intermediate intensity. We find however, that motif T2 exhibits greater regularity than motif T1 over the range $0.001 \leq \delta < 0.16$; before the optimal point B in Fig. 4a which occurs at $\delta_1=0.16$. That is, when $0.001 \leq \delta < 0.16$, the $R$ of the T2 motif is less than the $\gamma$ of the T1 motif (Fig. 4a). As $\delta_1$ approaches the optimal point $\delta_1=0.16$, the difference between the two curves becomes increasingly small. Interestingly, the optimal noise intensity $\delta_1^*$ and the corresponding absolute minimum in $R$, denoted $R^*$, are both the same for motifs T1 and T2. For example, in Fig.  4a, $\delta_1^*=0.162$ and $R^*=0.116$. Moreover, as $\delta_2$ increases further, the blue and black curves in Fig. 4a change at the same rate, and therefore the output regularity of motifs T1 and T2 is the same over this region.

The density functions of the ISIs presented in panels b and c of Fig.  4 corroborate our findings in Fig. 4a. We consider the density functions of the ISIs for three disparate intensities of the driving noise which correspond to points A, B, and C in Fig.  4a. Point A is a low level of noise with intensity $\delta_1=0.0023$ and corresponds to the orange curve, point B is a moderate level of noise with intensity $\delta_1=0.16$ and corresponds to the purple curve, and point C is a high level of noise with intensity $\delta_1=1.35$ and corresponds to the green curve. First, we find that the peak of the purple density function, which corresponds to the optimal noise intensity, is more pronounced in both panels b and c than those of the orange and green density functions which correspond to noise intensities that are either too weak or too strong. This implies that the ISIs at point B are more tightly grouped, and agrees with our results in Fig. 4a. Furthermore, the orange density function in panel c displays a much more pronounced peak (i.e. larger peak and smaller half-width) than its counterpart, the orange density function in panel b, which has relatively broad bimodal peaks. This suggests that the output of the T2 motif is more regular than that of the T1 motif over the low intensity range of $\delta_1$. Conversely, the purple and green density functions appear to be equivalent across panel b and c of Fig. 4. This is consistent with Fig.  4a, which suggests that sufficiently strong intensities of the driving noise can eliminate the differences in the output regularity of the T1 and T2 motifs. 

Overall, we observe that the regularity of the output of both motifs T1 and T2 increases as the noise intensity increases from low levels until it reaches an optimal point. Before this optimal point, the T2 motif shows greater output reliability (as measured by a lower value of $R$). And Beyond this optimal point, the regularity of the T1 and T2 motifs is identical, and with the increment of $\delta_1$, the regularity begins to decline as the noise intensity becomes overpowering.

\subsection{The effects of network connectivity on network synchrony and output regularity}
\begin{figure}
    \centering
    \includegraphics[scale=0.7]{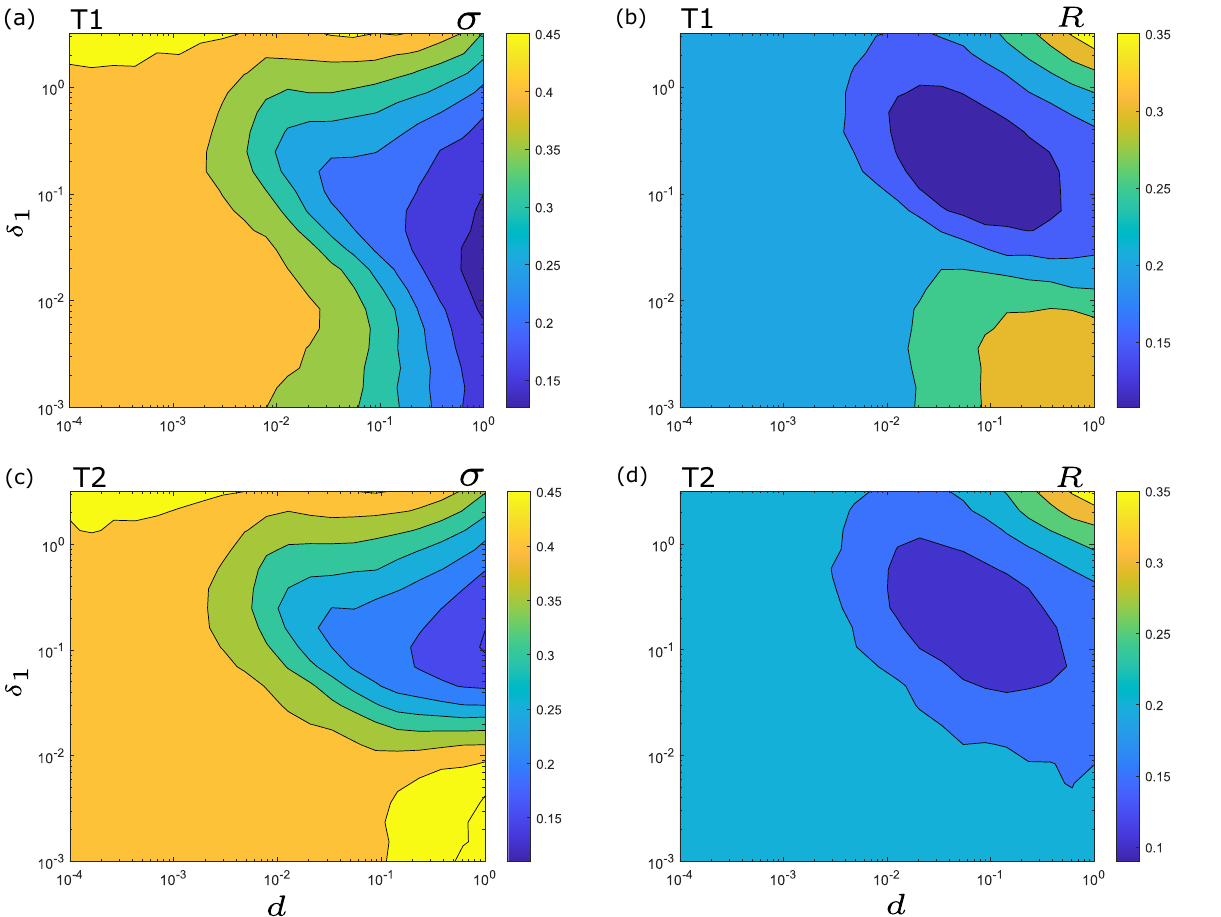}
    \caption{Contour maps of $\sigma$ and $R$ vs. noise intensity $\delta_1$ vs. coupling strength, $d$. (a): $\sigma$ vs. $\delta_1$ vs. $d$ for motif T1.(b): $R$ vs. $\delta_1$ vs. $d$ for motif T1. (c): $\sigma$ vs. $\delta_1$ vs. $d$ for motif T2. (d): $R$ vs. $\delta_1$ vs. $d$ for motif T2. In (a) and (b) $d=d_{3,1}=d_{3,2}=d_{2,1}$. In (c) and (d) $d=d_{3,1}=-d_{3,2}=d_{2,1}$. Other parameters are $\alpha=-0.2$, $\gamma=-0.2$, $\omega_0=2$, $\omega_1=0$, $\delta_2=\delta_3=0.01$, and $\lambda_0=-0.1$.}
    \label{fig:coupling}
\end{figure}
In Secs. 3.2 and 3.3, we investigated the impact of the driving noise intensity, $\delta_1$, on network synchronization and output regularity with fixed coupling strengths $d_{3,1}$, $d_{3,2}$, and $d_{2,1}$. Next we will consider the effects of the network connectivity on network synchronization and output regularity by considering a variable coupling strength $d$, where $d=d_{3,1},d_{3,2},d_{2,1}$ for the T1 motif, and $d=d_{3,1},-d_{3,2},d_{2,1}$ for the T2 motif. To quantify network synchronization we use the root mean square deviation, $\sigma$, in Eq. \ref{sigma} and to quantify output regularity we use the coefficient of variation, $R$, in Eq. \ref{coeff variation}. Our findings are presented in Fig. 5, in the form of contour maps of $\sigma$ and $R$ as functions of $\delta_1$ and $d$ for the T1 motif (panels a and b, respectively) and T2 motif (panels c and d, respectively). As indicated by the colour bars, warmer colours correspond to larger values of both $\sigma$ and $R$.

First, we see that there exist threshold values of the common coupling strength, $d$, denoted $d_T$, such that when $d \leq d_T$, changes in $\delta_1$ have no significant effect on network synchronization ($\sigma$) and output regularity ($R$). In the case of network synchronization, we find that $d_T \approx 0.02$ for both T1 and T2 motifs (Fig. 5a and 5c) and for output regularity, $d_T \approx 0.03$ for both the T1 and T2 motifs (Fig. 5b and 5d). Furthermore, we see that the contour maps in panels c and d, which correspond to motif T2, are notably different than the contour maps in panels a and b, which correspond to motif T1, when the noise intensity is weak (e.g. $\delta_1<0.01$ for $\sigma$ and $\delta_1<0.03$ for $R$ in Fig. 5) and the coupling strength is strong (e.g. $d>0.01$ for $\sigma$ and $d>0.02$ for $R$); or simply, the bottom right corners of the contour maps in Fig. 5. We find that the T1 motif has a smaller $\sigma$ in this region relative to the T2 motif, but a larger $R$. In other words, the T1 motif shows a greater propensity for network synchrony but a lesser propensity for output regularity than the T2 motif in this region---clearly, the converse holds as well. Moreover, we see that the differences in $\sigma$ and $R$ disappear as $\delta_1$ increases, which is consistent with our findings in the sections (upper right region of Fig. 5).
\begin{figure}
    \centering
   \includegraphics[scale=0.65]{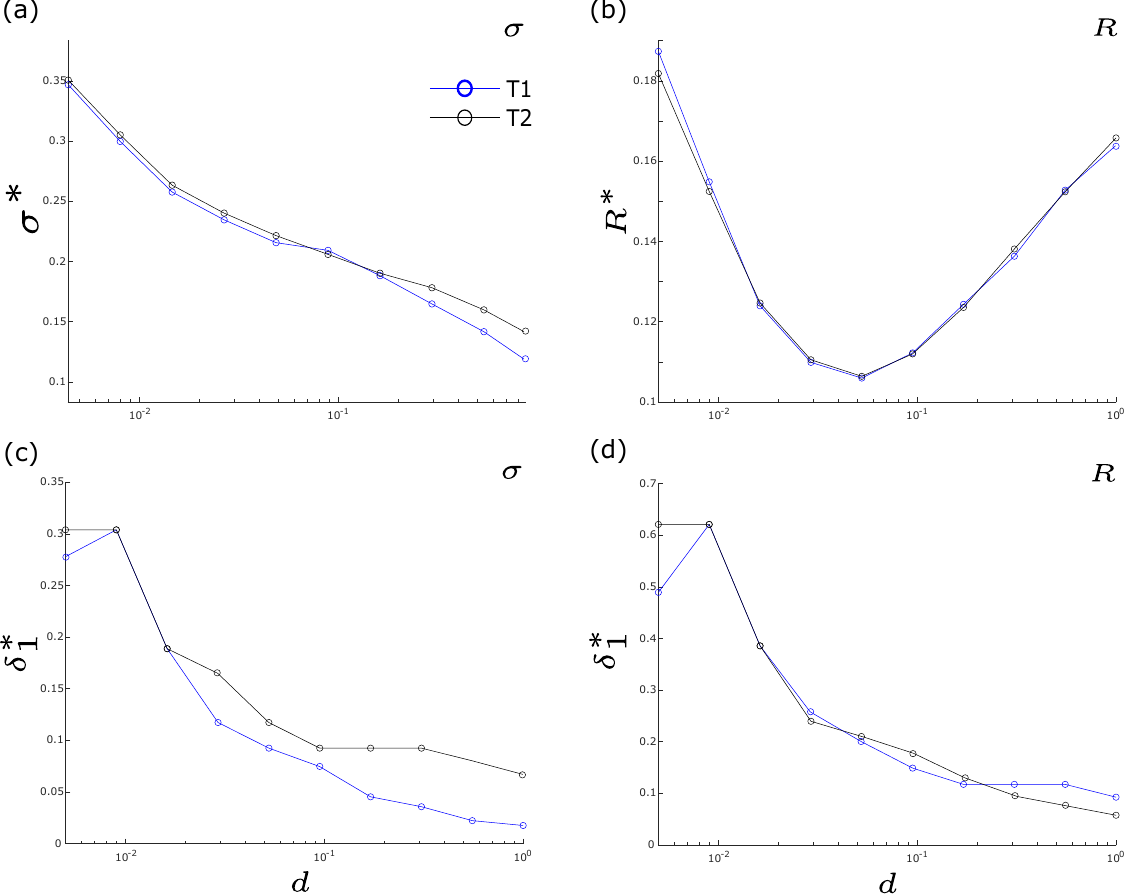}
   \caption{Optimal $\sigma$ and $R$, and their optimal noise intensities vs. network connectivity. Panels (a) and (c) show $\sigma^*$ and corresponding $\delta_1^*$ vs. $d$, respectively, while panels (b) and (d) show $R^*$ and corresponding $\delta_1^*$ vs. $d$, respectively. Blue and black curves correspond to T1 and T2 motifs, respectively. Parameters are: $\alpha=-0.2$, $\gamma=-0.2$, $\omega_0=2$, $\omega_1=0$, $\delta_2=\delta_3=0.01$, $\lambda_0=-0.1$, $d=d_{3,1}=d_{3,2}=d_{2,1}$ for the T1 motif, and $d=d_{3,1}=-d_{3,2}=d_{2,1}$ for the T2 motif.}
     \label{fig:optimal}
 \end{figure}
 
The contours in Fig. 5 suggest that the noise intensity required for both motifs to optimize both network synchrony ($\sigma$) and/or output regularity ($R$) is dependent on $d$. To consider the effects the network connectivity, $d$, more comprehensively, we compute $\sigma^*$ and $R^*$ and their corresponding $\delta_1^*$ values as functions of $d$. Our findings are presented in Fig. 6, where: panel (a) displays $\sigma^*$ vs. $d$; (b) displays $R^*$ vs. $d$; and (c) and (d) display their corresponding optimal noise intensities, $\delta_1^*$, vs. $d$, respectively. For motifs T1 and T2, we find that $\sigma^*$ is negatively correlated with $d$. This suggests that the optimal level of network synchrony increases as the network connectivity, $d$, increases. Additionally, in Fig. 6a, the $\sigma^*$ curve for the T1 motif is less than or equal to the $\sigma^*$ curve for the T2 motif. This is consistent with the results presented in Sec. 3.2 which highlight the T1 motif's propensity for greater network synchrony. Furthermore, Fig. 6c reveals that the optimal noise intensity, $\delta_1^*$ is negatively correlated with $d$. That is, as network connectivity increases, the intensity of noise required to optimize network synchrony decreases. In addition, the $\delta_1^*$ curve for motif T1 in Fig. 6c (blue curve) is less than or equal to the $\delta_1^*$ curve for the motif T2 (black curve). This suggests that the T1 motif requires weaker noise to maximize network synchronization than to the T2 motif (as in Fig. 3a).

From Fig. 6b we see the $R^*$ is a non-monotone function of $d$. For both motifs T1 (black curve) and T2 (blue curve), $R^*$ is decreasing when $d<0.05$, reaches a minimum at $d=0.05$, and then increases for $d>0.05$. When $d=0.05$, both motifs have the same optimal output regularity: $R^*=0.108$. Indeed, there is no clear difference between the $R^*$ curves for the T1 and T2 motifs in Fig. 6b. This corroborates the findings in Fig. 4, which show that the T1 and T2 motifs have similar levels of optimal output regularity. Fig. 6d, which displays the optimal noise intensities associated with the curves in Fig. 6b suggests $\delta_1^*$ decreases as $d$ increases for both types of motifs, and further, that there is no consistent difference between the $\delta_1$ required to minimize $R$ for the T1 motif and T2 motif. This is consistent with our results in Fig. 4.
 \subsection{The effects of $\lambda_0$ on network synchrony and output regularity}
  \begin{figure}
    \centering
   \includegraphics[scale=0.6]{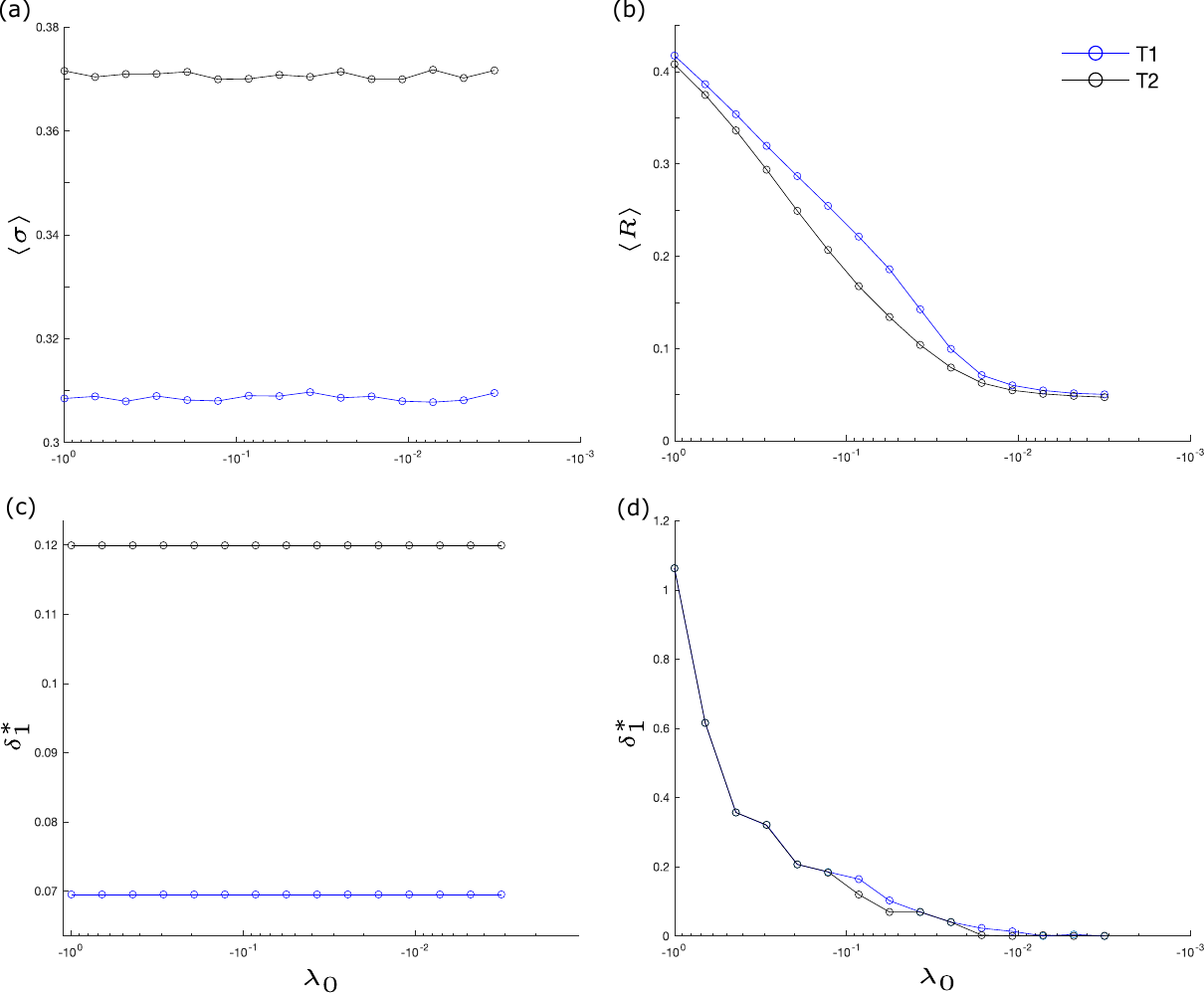}
   \caption{ (a) mean network synchronization,$\langle \sigma \rangle$, and (b) output reliability, $\langle R \rangle$, as functions of $\lambda_0$. (c) and (d) show the optimal noise intensities needed to optimize network synchronization and output reliability, respectively. Blue and black curves represent the T1 and T2 motifs, respectively. Other parameters are: $\alpha=-0.2$, $\gamma=-0.2$, $\omega_0=2$, $\omega_1=0$, $\delta_2=\delta_3=0.01$, $d_{3,1}=d_{3,2}=d_{2,1}=0.1$ for the T1 motif, and $d_{3,1}=-d_{3,2}=d_{2,1}=0.1$ for the T2 motif.}
     \label{fig:lambda}
 \end{figure}
 In the preceding sections, we have only considered $\lambda_0=-0.1$. However, many studies (e.g., \cite{yu2009constructive, yu2021noise}) show network synchronization in excitable systems is dependent on the distance of the critical parameter from the excitation threshold. Furthermore, we consider the effects of $\lambda_0$ on network synchronization and output regularity. To quantify network synchronization and output regularity we use the measures $\sigma$ in Eq. \ref{sigma} and $R$ in Eq. \ref{coeff variation}. We calculate $\sigma$ and $R$ within the range of $-1 \leq \lambda_0 \leq -0.001$ and $0.001\leq \delta_1 \leq 5$, and then compute the noise-averaged values of $\sigma$ and $R$, denoted by $\langle \sigma \rangle$ and $\langle R \rangle$, respectively, which allows us to consider $\langle \sigma \rangle$ and $\langle R \rangle$ as functions of $\lambda_0$. Our findings are presented in Figs. 7a and 7b, where the blue curves in both figures represent motif T1 and the black curves represent motif T2. Surprisingly, we find that there is no correlation between $\langle \sigma \rangle$ and $\lambda_0$ (Fig. 7a), both $\langle \sigma \rangle$ curves in Fig. 7a are constant with $\langle \sigma \rangle \approx 0.308$ for motif T1 and $\langle \sigma \rangle \approx 0.372$ for motif T2. That is, network synchronization is not correlated with the distance of $\lambda_0$ from the excitation threshold. Our results further reveal that $\langle \sigma \rangle$ is lower for the T1 motif relative to the T2 motif across the range $-1 \leq \lambda_0 \leq -0.08$. This finding is in agreement with our previous findings that suggest that motif T1 exhibits a greater degree of network synchronization than motif T2.

On the other hand, from Fig. 7b we see that there is a correlation between $\langle R \rangle$ and $\lambda_0$. For both motifs we find that $\langle R \rangle$ decreases linearly as $\lambda_0$ approaches the excitation threshold and then relaxes toward a minimum when $-0.01 \leq \lambda_0 < 0$, after which moving the networks closer to the excitation threshold does not affect $\langle R \rangle$. Our results further indicate that the black curve in Fig. 7b which corresponds to the T2 motif is always below the blue curve in Fig. 7b which corresponds to the T1 motif. This is in agreement with our previous findings, which suggest that the T2 motif exhibits greater output regularity than the T1 motif over the weak intensity range (e.g. Fig. 4a).\par

Next we consider how $\lambda_0$ affects the optimal driving noise intensity, $\delta_1^*$, needed to minimize $\sigma$ and $R$. Figs. 2c and 2d display the optimal driving noise intensities which correspond to $\sigma^*$ and $R^*$, respectively, as functions of $\lambda_0$. The T1 motifs are represented by blue curves and the T2 motifs are represented by black curves. We find that the $\delta_1^*$ corresponding $\sigma$ is a constant function of $\lambda_0$ for both the T1 and T2 motifs. Namely, $\delta_1^* \approx 0.07$ for the T1 motif and $\approx 0.12$ for the T2 motif (Fig. 7c). These results are consistent with our earlier findings in Fig. 3a, which suggest that the T1 motif requires a lesser intensity of $\delta_1$ to minimize $\sigma$ relative to the T2 motif. Conversely, in Fig. 7d we find that the $\delta_1$ needed to minimize $R$ is a decreasing function of $\lambda_0$, and as in Fig. 7b, the rate of change of the $\delta_1^*$ decreases as $\lambda_0 $ approaches the excitation threshold. Finally, the noise intensities required to optimize output regularity are approximately the same for both the T1 and T2 type motifs.

\section{Discussion}
We explore the impact of noise on network synchronization and output regularity within the context of three-neuron FFL motifs. Our investigation focuses exclusively on two distinct motif types: T1, which has purely excitatory connections (as depicted in Fig. 1a); and T2, which has a combination of both excitatory and inhibitory connections (as illustrated in Fig. 1b). We choose to analyze these specific motifs due to their prevalence in neuronal networks \cite{reigl2004search, Song2005}. Nevertheless, there are other three-neuron motifs that are also prevalent within the brain. More generally, neural motifs (i.e., frequently recurring wiring patterns) are thought to serve important functions within the brain, and may be viewed as the fundamental building blocks of larger, more complex neuronal networks \cite{Song2005, alon2007network, sporns2004motifs}.

We model the dynamics of each oscillator using the $\lambda-\omega$ system, the canonical model for the normal form of a HB, since is a well-established framework that captures the dynamics of the transition from quiescence to oscillations observed in real neurons \cite{izhikevich2000neural}. We consider the excitable regime, which is quiescent in the absence of noise but can be excited by the addition of an intrinsic noise. Our result show that the addition of a pure noise stimulus can maximize both network synchronization and output reliability at an intermediate intensity. Our findings agree with existing studies which find CR and other noise-induced effects in similar network motifs \cite{lu2019,guo2009,lou2014stochastic}.

Furthermore, our results indicate that the T1 motif shows a greater propensity for noise-induced synchronization; it exhibits a greater degree of synchronicity than the T2 motif and requires a comparatively lower noise intensity to optimize network synchrony. Conversely, we found that the T2 motif displays greater output reliability over the weak to intermediate range than the T1 motif but does not require a lesser intensity of noise to reach an optimal point. Importantly, our observations indicate that these differences between the two motifs are robust in that they persist over a wide range of noise intensities and parameter regimes, and only vanish when the intrinsic noise stimulus becomes excessively strong or the network connectivity becomes excessively weak. These results may in part be understood through our results in Fig. 2, which indicate: (i) the T2 motif shows both in-phase and anti-phase synchronization, whereas the T1 motif shows only in-phase synchronization (Fig. 2a, 2d); and (ii) the amplitude of the T2 oscillator is larger and appears to be more stable than the amplitude of the output oscillator of the T1 motif (Fig. 2c, 2f). Since we consider in-phase synchronization it is clear why (i) may lead motif T1 to exhibit a greater propensity for network synchrony than the T2 motif. On the other hand, from (ii) we find that the oscillations of $x_3$ of motif T2 have more pronounced and larger peaks than those which correspond to motif T1. This makes it easier to distinguish between spurious subthreshold oscillations and actual peaks and may therefore promote greater output regularity (or reliability) as we consider regularity in the context of the ISI variability.

Our results further suggest that network connectivity is positively correlated with synchronization and negatively correlated with the intensity of noise required to maximize synchronization, which is in agreement with existing studies (e.g., \cite{lu2019,guo2009,ge2020propagation}). Conversely, we find that the output regularity reaches an optimum point and an intermediate level of network connectivity ($d=0.05$; Fig. 6b). Finally, we find that moving our model closer to the excitation threshold has no effect on network synchronization or the amount of noise needed to maximize network synchrony, whereas it can indeed enhance output regularity and decrease the level of noise needed to optimize output regularity.

Overall, we find that the T1 (excitatory FFL) and T2 (excitatory and inhibitory FFL) motifs, which are common in biological networks, differ in terms of their respective propensities for network synchronization and output regularity. Our results emphasize the functional importance of neural motifs and the diverse roles that they may hold within the brain. Possible extensions of our work could include: (a) an investigation into how uneven interactions (like coupling and noise) impact emergent dynamics in neural motifs; and (b) an analysis of the functional differences between neural motifs embedded within larger networks; in addition to an examination of the stochastic dynamics networks composed of such motifs.

\bibliographystyle{elsarticle-num}
\bibliography{bib}

\begin{thebibliography}{10}
\expandafter\ifx\csname url\endcsname\relax
  \def\url#1{\texttt{#1}}\fi
\expandafter\ifx\csname urlprefix\endcsname\relax\def\urlprefix{URL }\fi
\expandafter\ifx\csname href\endcsname\relax
  \def\href#1#2{#2} \def\path#1{#1}\fi

\bibitem{Mangan2003}
S.~Mangan, U.~Alon, Structure and function of the feed-forward loop network
  motif, Proceedings of the National Academy of Sciences of the United States
  of America 100~(21) (2003) 11980--11985.
\newblock \href {https://doi.org/10.1073/pnas.2133841100}
  {\path{doi:10.1073/pnas.2133841100}}.

\bibitem{Shen2002}
S.~Shen-Orr, R.~Milo, S.~Mangan, U.~Alon, Network motifs in the transcriptional
  regulation network of escherichia coli, Nature Genetics 31~(1) (2002) 64--68.
\newblock \href {https://doi.org/10.1038/ng881} {\path{doi:10.1038/ng881}}.

\bibitem{reigl2004search}
M.~Reigl, U.~Alon, D.~B. Chklovskii, Search for computational modules in the c.
  elegans brain, BMC biology 2 (2004) 1--12.

\bibitem{Song2005}
S.~Song, P.~J. Sj{\"o}str{\"o}m, M.~Reigl, S.~Nelson, D.~B. Chklovskii, Highly
  nonrandom features of synaptic connectivity in local cortical circuits, PLoS
  biology 3~(3) (2005) e68.

\bibitem{Milo2002}
R.~Milo, S.~Shen-Orr, S.~Itzkovitz, N.~Kashtan, D.~Chklovskii, U.~Alon, Network
  motifs: Simple building blocks of complex networks, Science 298~(5594) (2002)
  824--827.
\newblock \href {https://doi.org/10.1126/science.298.5594.824}
  {\path{doi:10.1126/science.298.5594.824}}.

\bibitem{Lee2002}
T.~I. Lee, N.~J. Rinaldi, F.~Robert, D.~T. Odom, Z.~Bar-Joseph, G.~K. Gerber,
  N.~M. Hannett, C.~T. Harbison, C.~M. Thompson, I.~Simon, Microarray analysis
  of complex biological processes, Science 298 (2002) 799–804.
\newblock \href {https://doi.org/10.1126/science.1072981}
  {\path{doi:10.1126/science.1072981}}.

\bibitem{Boyer2005}
L.~A. Boyer, T.~I. Lee, M.~F. Cole, S.~E. Johnstone, S.~S. Levine, J.~P.
  Zucker, R.~A. Young, Core transcriptional regulatory circuitry in human
  embryonic stem cells, Cell 122~(6) (2005) 947–956.
\newblock \href {https://doi.org/10.1016/j.cell.2005.08.020}
  {\path{doi:10.1016/j.cell.2005.08.020}}.

\bibitem{tsang2007}
J.~Tsang, J.~Zhu, A.~van Oudenaarden, Microrna-mediated feedback and
  feedforward loops are recurrent network motifs in mammals, Molecular Cell
  26~(5) (2007) 753–767.
\newblock \href {https://doi.org/10.1016/j.molcel.2007.05.018}
  {\path{doi:10.1016/j.molcel.2007.05.018}}.

\bibitem{sporns2004}
O.~Sporns, R.~Kötter, Motifs in brain networks, PLoS Biol 2~(11) (2004) e369.
\newblock \href {https://doi.org/10.1371/journal.pbio.0020369}
  {\path{doi:10.1371/journal.pbio.0020369}}.

\bibitem{milo2002network}
R.~Milo, S.~Shen-Orr, S.~Itzkovitz, N.~Kashtan, D.~Chklovskii, U.~Alon, Network
  motifs: simple building blocks of complex networks, Science 298~(5594) (2002)
  824--827.

\bibitem{Eggermont1998}
J.~Eggermont, Representation of spectral and temporal sound features in three
  cortical fields of the cat, J Neurophysiol 80 (1998) 2743–2764.
\newblock \href {https://doi.org/10.1152/jn.1998.80.6.2743}
  {\path{doi:10.1152/jn.1998.80.6.2743}}.

\bibitem{Middleton2006}
J.~W. Middleton, A.~Longtin, J.~Benda, L.~Maler, The cellular basis for
  parallel neural transmission of a high-frequency stimulus and its
  low-frequency envelope, Proceedings of the National Academy of Sciences
  103~(39) (2006) 14596–14601.
\newblock \href {https://doi.org/10.1073/pnas.0604103103}
  {\path{doi:10.1073/pnas.0604103103}}.

\bibitem{Gui2016}
R.~Gui, Q.~Liu, Y.~Yao, Noise decomposition principle in a coherent
  feed-forward transcriptional regulatory loop, Front Physiol 7 (2016) 600.
\newblock \href {https://doi.org/10.3389/fphys.2016.00600}
  {\path{doi:10.3389/fphys.2016.00600}}.

\bibitem{Krauss2019}
P.~Krauss, K.~Prebeck, A.~Schilling, C.~Metzner, Recurrence resonance in
  three-neuron motifs, Front Comput Neurosci 13 (2019) 64.
\newblock \href {https://doi.org/10.3389/fncom.2019.00064}
  {\path{doi:10.3389/fncom.2019.00064}}.

\bibitem{wang2009}
Q.~Wang, M.~Perc, Z.~Duan, G.~Chen, Synchronization transitions on scale-free
  neuronal networks due to finite information transmission delays, Physical
  Review E 80~(2) (2009) 026206.

\bibitem{gao2001}
Z.~Gao, B.~Hu, G.~Hu, Stochastic resonance of small-world networks, Physical
  Review E 65~(1) (2001) 016209.

\bibitem{rosenblum2001phase}
M.~Rosenblum, A.~Pikovsky, J.~Kurths, C.~Sch{\"a}fer, P.~A. Tass, Phase
  synchronization: from theory to data analysis, in: Handbook of biological
  physics, Vol.~4, Elsevier, 2001, pp. 279--321.

\bibitem{gang1993stochastic}
H.~Gang, T.~Ditzinger, C.-Z. Ning, H.~Haken, Stochastic resonance without
  external periodic force, Physical review letters 71~(6) (1993) 807.

\bibitem{pikovsky1997coherence}
A.~S. Pikovsky, J.~Kurths, Coherence resonance in a noise-driven excitable
  system, Physical Review Letters 78~(5) (1997) 775.

\bibitem{mormann2000mean}
F.~Mormann, K.~Lehnertz, P.~David, C.~E. Elger, Mean phase coherence as a
  measure for phase synchronization and its application to the eeg of epilepsy
  patients, Physica D: Nonlinear Phenomena 144~(3-4) (2000) 358--369.

\bibitem{bonsel2022control}
F.~B{\"o}nsel, P.~Krauss, C.~Metzner, M.~E. Yamakou, Control of noise-induced
  coherent oscillations in three-neuron motifs, Cognitive Neurodynamics 16~(4)
  (2022) 941--960.

\bibitem{lu2019phase}
L.~Lu, C.~Bao, M.~Ge, Y.~Xu, L.~Yang, X.~Zhan, Y.~Jia, Phase noise-induced
  coherence resonance in three dimension memristive hindmarsh-rose neuron
  model, The European Physical Journal Special Topics 228 (2019) 2101--2110.

\bibitem{yu2009constructive}
G.~Yu, M.~Yi, Y.~Jia, J.~Tang, A constructive role of internal noise on
  coherence resonance induced by external noise in a calcium oscillation
  system, Chaos, Solitons \& Fractals 41~(1) (2009) 273--283.

\bibitem{yu2021noise}
N.~Yu, G.~Jagdev, M.~Morgovsky, Noise-induced network bursts and coherence in a
  calcium-mediated neural network, Heliyon (2021) e08612.

\bibitem{alon2007network}
U.~Alon, Network motifs: theory and experimental approaches, Nature Reviews
  Genetics 8~(6) (2007) 450--461.

\bibitem{sporns2004motifs}
O.~Sporns, R.~K{\"o}tter, Motifs in brain networks, PLoS biology 2~(11) (2004)
  e369.

\bibitem{izhikevich2000neural}
E.~M. Izhikevich, Neural excitability, spiking and bursting, International
  journal of bifurcation and chaos 10~(06) (2000) 1171--1266.

\bibitem{lu2019}
L.~Lu, Y.~Jia, J.~B. Kirunda, Y.~Xu, M.~Ge, Q.~Pei, L.~Yang, Effects of noise
  and synaptic weight on propagation of subthreshold excitatory postsynaptic
  current signal in a feed-forward neural network, Nonlinear Dynamics 95~(2)
  (2019) 1673--1686.

\bibitem{guo2009}
D.~Guo, C.~Li, Stochastic and coherence resonance in feed-forward-loop neuronal
  network motifs, Physical Review E 79~(5) (2009) 051921.

\bibitem{lou2014stochastic}
X.~Lou, Stochastic resonance in neuronal network motifs with ornstein-uhlenbeck
  colored noise, Mathematical Problems in Engineering 2014 (2014).

\bibitem{ge2020propagation}
M.~Ge, Y.~Jia, L.~Lu, Y.~Xu, H.~Wang, Y.~Zhao, Propagation characteristics of
  weak signal in feedforward izhikevich neural networks, Nonlinear Dynamics 99
  (2020) 2355--2367.

\end{thebibliography}
\end{document}